\documentclass[12pt]{article}

      \usepackage{latexsym}
         \usepackage[reqno, namelimits, sumlimits]{amsmath}
         \usepackage{amssymb, amsfonts}
         \usepackage{amsthm}

 \newcommand{\beqn}{\begin{eqnarray}}
 \newcommand{\eeqn}{\end{eqnarray}}
 \newcommand{\be}{\begin{equation}}
 \newcommand{\ee}{\end{equation}}
 \newcommand{\ba}{\begin{array}}
 \newcommand{\ea}{\end{array}}
 \newcommand{\pa}{\partial}
 
 \newcommand{\ci}{\cite}
 \newcommand{\la}{\label}

\newcommand{\Om}{\Omega}

\newcommand{\ga}{\gamma}

\newcommand{\na}{\nabla}
\newcommand{\om}{\omega}

 \newcommand{\De}{\Delta}

\newcommand{\M}{\mathcal{M}}

\def\R{{\rm I\kern-.1567em R}}
\def\M{{\rm I\kern-.1567em M}}

\def\div {{\rm div\,}}

\def\dist{{\rm dist}}

 \newtheorem{theorem}{Theorem}[section]

 \newtheorem{lemma}[theorem]{Lemma}
 
 \newtheorem{remark}[theorem]{Remark}

 \newtheorem{pro}[theorem]{Proposition}

\title{H\"{o}lder Continuity of Solutions of 2D Navier-Stokes Equations with Singular Forcing}

 \author{Peter Constantin\footnote{The University of Chicago} \and
  Gregory Seregin\footnote{Oxford University}}

\begin{document}
\maketitle

\vspace{1cm} \centerline{Dedicated to Nina Nikolaevna Uraltseva}
  \vspace{1cm}
 \noindent
 {\bf Abstract } We discuss the regularity of solutions of 2D incompressible Navier-Stokes
 equations forced by singular forces. The problem is motivated by the study of complex
 fluids modeled by the Navier-Stokes equations coupled to
a nonlinear Fokker-Planck equation describing microscopic corpora embedded in the fluid.
This leads naturally to bounded added stress and hence to
$W^{-1,\infty}$ forcing of the Navier-Stokes equations.

 \vspace {1cm}

\noindent {\bf 1991 Mathematical subject classification (Amer. Math.
Soc.)}: 35K, 76D.


\setcounter{equation}{0}
\section{Introduction  } We discuss the regularity of solutions of 2D incompressible Navier-Stokes
 equations forced by singular forces. The problem is motivated by the study of complex
 fluids modeled by the Navier-Stokes equations coupled to
a nonlinear Fokker-Planck equation describing microscopic corpora
embedded in the fluid. This leads naturally to bounded added stress
and hence to $W^{-1,\infty}$ forcing of the Navier-Stokes equations.
A more detailed description of the problem in question, together with an
application of the results in the present paper can be found in our
forthcoming paper \ci{CoSe}.

In this paper we focus on the 2D Navier-Stokes issues. The
global existence of energy solutions and their uniqueness are well
known as classical results of J. Leray for the Cauchy problem and O.
Ladyzhenskaya for initial boundary value problems in bounded
domains. These results remain to be true  for singular forces as
well.

The regularity of energy solutions with relatively smooth forces is
also known. Regularity can be established, for instance,  by scalar multiplication of the Navier-Stokes equation by the Stokes operator of the velocity field,
integration by parts, and application of Ladyzhenskay's
inequality
$$\|u\|^2_{L^4(\mathbb R^2)}\leq
\sqrt{2}\|u\|_{L^2(\mathbb R^2)}\|\na u\|_{L^2(\mathbb R^2)},\qquad
\forall u\in C^\infty_0(\mathbb R^2).$$ This
procedure yields summability of the second spatial derivatives. Further regularity can be obtained perturbatively, with the help of the linear theory.

The regularity of energy solutions with singular forcing is limited. The best one can expect is H\"older continuity of the velocity field. We prove H\"{o}lder continuity at a local level, in both space and in time. We  assume that our local solution has finite energy and the pressure field is in
$L^2$. This latter assumption seems restrictive: we are not able to
justify it for general initial boundary value
problems with reasonable singular forcing. The assumption is however satisfied
in the absence of boundaries, i.e., for the Cauchy problem in the whole space and for the initial value problem on
the torus. We briefly explain in this paper how the local regularity results
can be applied to the Cauchy problem in the whole space.

In our proof, the H\"older continuity
of the velocity field depends quantitatively on the modulus of
continuity of the function $\om \mapsto \int\limits_\om |u|^4dz$.
In order to be able to apply this regularity result to coupled systems or to families of Navier-Stokes systems, this modulus of continuity needs to be a priori uniformly controlled.
We achieve this in the absence of boundaries by obtaining higher integrability
of the velocity, $u\in L^{\infty}(dt; L^r(dx))$, $r\ge 4$.
In order to obtain the higher integrability we prove the generalized
Ladyzhenskaya inequality that reads
$$ \|u\|^2_{L^{2r}(\mathbb R^2)}\leq
\frac r{\sqrt{2}}\|u\|_{L^r(\mathbb R^2)}\|\na u\|_{L^2(\mathbb R^2)},\qquad \forall u\in C^\infty_0(\mathbb R^2)$$ for $r\geq 2$.
The proof is elementary and can be found in the Appendix.

\setcounter{equation}{0}
\section{Notation and Local Regularity Result}
 We assume that $\Om$ and $\Om_1$ are domains in
 $\mathbb R^2$ such that $\Om_1\Subset\Om$ and $0<T_1<T$,
 and let $$Q=\Om\times (-T,0),\qquad Q_1=\Om_1\times (-T_1,0).$$
Parabolic balls will be denoted as $Q(z_0, R)=B(x_0,R)\times (t_0-R^2,t_0)$,
where $z_0=(x_0,t_0)$, $x_0\in\mathbb R^2$, $t_0\in\mathbb R$, and $B(x_0,R)$ is an open disk
in $\mathbb R^2$ having  radius $R$ and centered at the point $x_0$.\\
We use the following notation for mean values:
$$(f)_{z_0,R}=\frac 1{|Q(z_0,R)|}\int\limits_{Q(z_0,R)}f(z)dz,\qquad
[p]_{x_0,R}=\frac 1{|B(x_0,R)|}\int\limits_{B(x_0,R)}p(x)dx.$$
$L^p(\Om)$ and $W^{l,p}(\Om)$ stand for usual Lebesgues and Sobolev spaces
of functions defined $\Om$, and the norm of the Lebesgues space is denoted
by $\|\cdot\|_{m,\Om}$. For the forcing we are going use a functional space
$M_{2,\ga}(Q)$ with parameter $0\leq \ga<1$ and seminorm
$$\|f\|_{M_{2,\ga}(Q)}=\sup\limits_{Q(z_0,R)\subset Q}R^{1-\ga}
\Big(\frac 1{|Q(z_0,R)|}\int\limits_{Q(z_0,R)}|f(z)-(f)_{z_0,R}|^2dz\Big)^\frac 12<\infty.$$
We denote by $c$ all positive universal constants. Our regularity result can be formulated as follows.
\begin{theorem}\la{nrt1} Assume that we are given functions
\be\la{nr1}u\in L^4(Q;\mathbb R^2),\quad p\in L^2(Q),\quad F \in M_{2,\ga}(Q;\mathbb M^{2\times 2})\ee
with $0\leq \ga<1$, satisfying the Navier-Stokes equations
\be\la{nr2}\pa_t u+u\cdot \na u-\De u +\na p = -\div F,\quad \div u=0\ee
in $Q$ in the sense of distributions.

Then
\be\la {nt3} u\in C^\ga(\overline{Q}_1)\ee
if $0<\ga<1$ and
\be\la {nr4} u \in BMO(Q_1)\ee
if $\ga=0$.
\end{theorem}
\begin{remark}\la{nrr2} The H\"older continuity and the BMO space are defined with
respect to parabolic metrics.\end{remark}
\begin{remark}\la{nrr3} The corresponding norms are estimated in terms of the quantities
$\|u\|_{4,Q}$, $\|p\|_{2,Q}$, $\|F\|_{M_{2,\ga}(Q)}$, $\dist(\Om_1,\pa\Om)$, $T-T_1$, and
the modulus of continuity of the function $\om \mapsto \int\limits_\om |u|^4dz$. \end{remark}
Several additional results can be proved by means of minor modifications
of the proof of Theorem \ref{nrt1}. Before stating one them, we define usual
energy spaces for the 2D Navier-Stokes equations. Let $H$ and $V$ be completions of the set of all divergence-free vector fields from $C^\infty_0(\mathbb R^2;\mathbb R^2)$ with respect to the $L^2$ norm and the Dirichlet integral, respectively.
\begin{pro}\la{nrp4} Let $u\in L^{\infty}(0,T;H)\cap L^2(0,T;V)$, $p\in L^2(0,T;L^2(\mathbb R^2))$ be a  solution of the Cauchy problem
\be\la{nr5}\pa_tu+u\cdot\na u-\De u+\na p=-\div F,\qquad \div
u=0,\ee \be\la{nr6}v(\cdot,0)=a(\cdot)\in H,\ee where $F\in
L^q(Q_T;\mathbb M^{2\times 2})\cap L^2(Q_T;\mathbb M^{2\times 2})$
with $q>4$ and $Q_T=\mathbb R^2\times (0,T)$.

Then, given $0<s \le T$, there exists a constant $C$ depending only
on $s$, the norms of $F$ in $L^q(Q_T; \mathbb M^{2\times 2})$ and in
$L^2(Q_T;\mathbb M^{2\times 2})$ and the norm of $a$ in $H$ such
that \be\la{np6} \|u\|_{L^{\infty}(\mathbb R^2\times (s,T))} \le
C.\ee Moreover, the function $u$ is H\"older continuous in $\mathbb
R^2\times [s,T]$ with exponent $\ga=1-\frac 4q$.\end{pro}
\begin{remark}\la{rnr5} The existence  and uniqueness of a solution to the Cauchy problems
(\ref{nr5}) and (\ref{nr6}) with above properties is well known, see \ci{La}.\end{remark}
\begin{remark}\la{rnr6} The same statement is valid in the case of periodic boundary conditions. More generally, it is true as long as the pressure field is in $L^2$. \end{remark}

\setcounter{equation}{0}
\section{Proof of Theorem \ref{nrt1} }
We are going to analyze differentiability properties of the velocity field $u$
in terms of the following functionals:
$$\Phi(u;z_0,\varrho)=\Big(\int\limits_{Q(z_0,\varrho)}
|u-(u)_{z_0,\varrho}|^4dz\Big)^\frac 12,\quad \Psi(u;z_0,\varrho)=\Big(\int\limits_{Q(z_0,\varrho)}
|u|^4dz\Big)^\frac 12,$$ $$D(p;z_0,\varrho)=\int\limits_{Q(z_0,\varrho)}|p-[p]_
{x_0,\varrho}|^2dz.$$

The following two statements are well-known.
\begin{lemma}\la{pl1} Let the function $v\in L^4(Q(z_0,R))$ satisfy the heat equation
$$\pa_tv-\De v=0$$ in $Q(z_0,R)$. Then
\be\la{p1}\Phi(v;z_0,\varrho)\leq c\Big(\frac \varrho R\Big)^4\Phi(v;z_0,R)\ee for all $0<\varrho\leq R$.\end{lemma}
\begin{lemma}\la{pl2} Given $G\in L^2(Q(z_0,R);\mathbb M^{2\times 2})$, there exists a unique function
$$w\in C([t_0-R^2,t_0];L^2(B(x_0,R);\mathbb R^2))\cap L^2([t_0-R^2,t_0];W^{1,2}(B(x_0,R);\mathbb R^2))$$
such that
$$\pa_tw-\De w=-\div G$$ in $Q(z_0,R)$
and $$w=0$$ on the parabolic boundary of $Q(z_0,R)$. Moreover, the function $w$ satisfies the estimates:
$$|w|^2_{2,Q(z_0,R)}\equiv \sup\limits_{t_0-R^2<t<t_0}\|w(\cdot,t)\|^2_{2,B(x_0,R)}+\|\na w\|^2_{2,Q(z_0,R)}$$\be\la{p2}\leq 2\|G\|^2_{2,Q(z_0,R)},\ee
\be\la{p3}\Phi(w;z_0,R)\leq c|w|^2_{2,Q(z_0,R)}.\ee\end{lemma}
The next couple of statements are about  some properties of the solutions of the system (\ref{nr2}).
\begin{lemma}\la{pl3}Under the assumptions of Theorem \ref{nrt1} we have
$$\Phi(u;z_0,\varrho)\leq c\Big\{\Big[\Big(\frac \varrho R\Big)^4+\Psi(u;z_0,R)\Big]\Phi(u;z_0,R)+$$
\be\la{p4}+D(p;z_0,R)+MR^{2+2\ga}\Big\}\ee
whenever $Q(z_0,R)\subset Q$ and $0<\varrho\leq R$. Here, $M=\|F\|^2_{M_{2,\ga}(Q)}$.\end{lemma}
\textsc{Proof}. Setting
$$G=F-(F)_{z_0,R}+(p-[p]_{x_0,R})\mathbb I+(u-(u)_{z_0,R})\otimes u$$
in Lemma \ref{pl2}, we get the following estimate for $w$
$$\Phi(w;z_0,R)\leq$$$$\leq c\int\limits_{Q(z_0,R)}\Big[|F-(F)_{z_0,R}|^2+
|p-[p]_{x_0,R}|^2+|u-(u)_{z_0,R}|^2|u|^2\Big]dz\leq$$
\be\la{p5}\leq c\Big[MR^{2+2\ga}+D(p;z_0,R)+\Psi(u;z_0,R)\Phi(u;z_0,R)\Big].\ee
Obviously, the function $v=u-w$ satisfies the heat equation. Then,  applying Lemma
\ref{pl1}, we find
$$\Phi(u-w;z_0,\varrho)\leq c\Big(\frac \varrho R\Big)^4\Phi(u-w;z_0,R).$$
The latter inequality gives us:
\be\la{p6}\Phi(u;z_0,\varrho)\leq c\Big[\Big(\frac \varrho R\Big)^4\Phi(u;z_0,R)
+\Phi(w;z_0,R)\Big].\ee
Combining (\ref{p5}) and (\ref{p6}), we arrive at (\ref{p4}) and thus Lemma \ref{pl3} is proved.
\begin{lemma}\la{pl4}Under the assumptions of Theorem \ref{nrt1}, we have the estimate
\be\la{p7}D(p;z_0,\varrho)\leq c\Big[\Big(\frac \varrho R\Big)^4D(p;z_0,R)+\Psi(u;z_0,R)\Phi(u;z_0,R)
+MR^{2+2\ga}\Big]\ee
whenever $Q(z_0,R)\subset Q$ and $0<\varrho\leq R$.\end{lemma}
\textsc{Proof}. The crucial part of the proof is the pressure decomposition
$$p=p_1+p_2,$$
where the first component $p_1$ satisfies
the identity
$$\int\limits_{B(x_0,R)}p_1\De \varphi dx=-\int\limits_{B(x_0,R)}\Big((u-(u)_{z_0,R})
\otimes u+F-(F)_{z_0,R}\Big):\na^2\varphi dx$$
where the test function $\varphi\in W^2_2(B(x_0,R))$ is subject to the Dirichlet boundary
condition: $\varphi=0$ on $\pa B(x_0,R)$. It is not difficult to show that such a function $p_1$ exists and obeys the estimate
$$\int\limits_{Q(z_0,R)}|p_1-[p_1]_{x_0,R}|^2dz\leq c\int\limits_{Q(z_0,R)}
|p_1|^2dz\leq $$\be\la{p8}\leq c\Big[\Psi(u;z_0,R)\Phi(u;z_0,R)
+MR^{2+2\ga}\Big].\ee The second counterpart of the pressure $p_2$ is a harmonic
function and thus satisfies the estimate
$$\int\limits_{B(x_0,\varrho)}|p_2-[p_2]_{x_0,\varrho}|^2dx\leq
c\Big(\frac \varrho R\Big)^4\int\limits_{B(x_0,R)}|p_2-[p_2]_{x_0,R}|^2dx$$
$$\leq c\Big(\frac \varrho R\Big)^4\int\limits_{B(x_0,R)}|p-[p]_{x_0,R}|^2dx+
c\int\limits_{B(x_0,R)}|p_1|^2dx$$
for any $0<\varrho\leq R$. Hence, by (\ref{p8}), we show
\be\la{p9}D(p_2;z_0,\varrho)\leq c\Big[\Big(\frac \varrho R\Big)^4D(p,z_0,R)+\Psi(u;z_0,R)\Phi(u;z_0,R)
+MR^{2+2\ga}\Big].\ee
Taking into account simple inequality
$$D(p;z_0,\varrho)\leq 2D(p_1;z_0,\varrho)+2 D(p_2;z_0,\varrho)\leq $$$$\leq \int\limits_{Q(z_0,R)}
|p_1|^2dz+2D(p_2;z_0,\varrho)$$
 we deduce the estimate (\ref{p7}) from (\ref{p8}) and (\ref{p9}).
 Lemma \ref{pl4} is proved.

Now, we pass to the  proof of Theorem \ref{nrt1}. Assuming that $Q(z_0,R)\subset Q$ and $\tau\in (0,1)$, we find
from (\ref{p4}) and (\ref{p7}) two inequalities:
$$\Phi(u;z_0,\tau^2R)\leq c(\tau^4+\Psi(u;z_0,\tau R))\Phi(u;z_0,\tau R)
+cD(p;z_0,\tau R)+$$$$+cM(\tau R)^{2+2\ga}$$
and
$$(1+c)D(p;z_0,\tau R)\leq (1+c)c\Big\{\tau^4D(p;bz_0,R)+\Psi(u;z_0,R)\Phi(u;z_0,R)
+$$$$+MR^{2+2\ga}\Big\}.$$
Adding the latter inequalities and introducing the new functional
$$\Theta(z_0,R)=\Phi(u;z_0,\tau R)+D(p; z_0,R),$$
we arrive at the basic estimate
$$\Theta(z_0,\tau R)\leq c(\tau^4+\Psi(u;z_0,\tau R))\Theta(z_0,R)+$$\be\la{p10}+
c\Psi(u;z_0,R)\Phi(u;z_0,R)+cM R^{2+2\ga}.\ee

It is not difficult to check validity of the following inequality
$$\Phi(u;z_0,\tau R)\leq c\Phi(u;z_0,R)$$
for any $\tau\in (0,1)$ and any $R>0$. Then from (\ref{p10}) it follows that
\be\la{p11}\Theta(z_0,\tau R)\leq c(\tau^4+\Psi(u;z_0, R))\Theta(z_0,R/\tau)+cM R^{2+2\ga}\ee
under assumption that $Q(z_0,R/\tau)\subset Q$.
Letting $\ga_1=(1+\ga)/2$ and choosing $\tau=\tau(\ga)\in (0,1)$ so that
\be\la{p12}c\tau^{2-2\ga_1}\leq 1/2,\ee
we can choose $R_0<\tau\min\{\dist(\Om_1,\pa\Om), \sqrt{T-T_1}\}$ so that
\be\la{p13}\Psi(u;z_0,R)<\tau^4\ee
for all $z_0\in Q_1$ and all $0<R\leq R_0$. It is here that the modulus of continuity of $\omega\mapsto \int_{\omega}|u|^4dxdt$ is used. So, summarizing all the above,
we have
\be\la{p14}\Theta(z_0,\tau R)\leq \tau^{2+2\ga_1}\Theta(z_0,R/\tau)+cM R^{2+2\ga}\ee
for all $z_0\in Q_1$ and all $0<R\leq R_0$. To reduce (\ref{p14}) to a known
iterative procedure, we let $\varrho=R/\tau$ and $\vartheta=\tau^2$. As a result, we find
\be\la{p15}\Theta(z_0,\vartheta \varrho)\leq \vartheta^{1+\ga_1}\Theta(z_0,\varrho)+cM \vartheta^{1+\ga}\varrho^{2+2\ga}\ee
for all $z_0\in Q_1$ and all $0<\rho\leq R_0/\tau$. The inequality (\ref{p15})
can be easily iterated, see \ci{Cam},
$$\Theta(z_0,\vartheta^k R_0/\tau)\leq \vartheta^{k(1+\ga)}\Big(\Theta(z_0, R_0/\tau)
+c_1MR_0^{2+2\ga}\Big)$$
for any $k\in \mathbb N$. Here and in what follows, all positive constants depending on $\ga$ only are denoted by $c_1$. The latter inequality implies
$$\Phi(u;z_0,\vartheta^k R_0)\leq \vartheta^{k(1+\ga)}\Big(\Theta(z_0, R_0/\tau)
+c_1MR_0^{2+2\ga}\Big) $$ for any $k\in \mathbb N$ and, hence,
\be\la{p16}\Phi(u;z_0,\varrho)\leq \varrho^{1+\ga}H\ee
for any $0<\varrho\leq R_0$, where
$$H=c_1\Big(\frac 1 {R_0}\Big)^{1+\ga}\Big[\Theta(z_0, R_0/\tau)
+MR_0^{2+2\ga}\Big].$$
Obviously, $H$ is a function of $R_0$, $\|u\|_{4,Q}$, $\|p\|_{2,Q}$, $M$, $\dist(\Om_1,\pa\Om)$, $T-T_1$, and $\ga$.

Now, our next step is to figure out how does $\Psi(u;z_0,R)$ depend on $R$. By (\ref{p16}),
we have
$$|(u)_{z_0,\varrho/2}-(u)_{z_0,\varrho}|\leq \frac c\varrho\Phi^{1/2}(u;z_0,\varrho)\leq cH^{1/2} \varrho^{-(1-\ga)/2}$$
for any $0<\varrho\leq R_0$. Therefore,
$$|(u)_{z_0,R_0/2^k}-(u)_{z_0,R_0}|\leq c H^{1/2}\sum\limits^{k-1}_{i=0}
\Big(\frac {R_0}{2^i}\Big)^{-(1-\ga)/2}=$$$$=cH^{1/2}\Big(\frac {R_0}{2^{k}}\Big)^{-(1-\ga)/2}\sum\limits^{k-1}_{i=0}\Big(\frac 1{2^{k-i}}\Big)^{(1-\ga)/2}\leq c_1H^{1/2}\Big(\frac {R_0}{2^{k}}\Big)^{-(1-\ga)/2}$$ for any $k\in \mathbb N$, or
\be\la{p17}|(u)_{z_0,\varrho}-(u)_{z_0,R_0}|\leq c_1 H^{1/2}\frac 1{\varrho^{(1-\ga)/2}}\ee for any $0<\varrho\leq R_0$. Proceeding
and making use of (\ref{p17}), we find
$$\Psi^{1/2}(u;z_0,\varrho)\leq \Phi^{1/2}(u;z_0,\varrho)+c\varrho
|(u)_{z_0,\varrho}|\leq c_1\varrho^{(1+\ga)/2}H^{1/2}+$$
$$+c\varrho|(u)_{z_0,R_0}|.$$
The latter implies
\be\la{p18}\Psi(u;z_0,\varrho)\leq \varrho^{(1+\ga)}H_1\ee
for any $0<\varrho\leq R_0$, where $H_1$ depends on the same arguments as $H$.

Coming back to the basic estimate (\ref{p10}) and taking into account (\ref{p12}), (\ref{p13}), (\ref{p16}),
and (\ref{p18}),
$$\Theta(z_0,\tau R)\leq \tau^{2+2\ga}\Theta(z_0,R)
+c(M+HH_1) R^{2+2\ga}$$
for any $0<R\leq R_0$. After iterations of the latter inequality, we show
$$\Theta(z_0,\tau^k R)\leq \tau^{k(2+2\ga)}\Big(\Theta(z_0,R_0)
+c(M+HH_1) R_0^{2+2\ga}\Big)$$
for any $k\in\mathbb N$. Consequently, we obtain:
\be\la{p19}\Phi(u;z_0,\varrho)\leq \varrho^{2+2\ga}H_2\ee
for any $z_0\in Q_1$ and for any $0<\varrho\leq R_0$ with $H_2$ depending on the same arguments
as $H$ and $H_1$. Finally, the H\"older continuity of $u$ subject to (\ref{p19}) follows from
known considerations, see, for example, \ci{DaPr} or \ci{GiaStr}. Theorem \ref{nrt1} is proved.

\noindent
\textsc{
Proof of Proposition \ref{nrp4}} Our first remark is that
\be\la{p20}L^q(Q_T)\subset M_{2,\ga}(Q_T)\ee if $\ga=1-\frac 4q$.
The second remark is that \be\la{p21}u\in L^4(Q_T)\ee and the
corresponding norm is bounded by $\|a\|_{2,\mathbb
R^2}+\|F\|_{2,Q_T}$.

 In order to handle the modulus of continuity  of
the function $\om \mapsto \int\limits_\om |u|^4dz$, we are going to
show that, for any $0<s<T$,
\be\la{p22}\sup\limits_{s<t<T}\|u(\cdot,t)\|_{4,\mathbb R^2}\leq
C(s,\|F\|_{4,Q_T},\|a\|_{2,\mathbb R^2}).\ee Assume that (\ref{p22})
has been already proved. Let the number $\tau\in (0,1)$ be defined
by (\ref{p12}). Then by (\ref{p22}) we can find a number
$0<R_1<\frac 12 \sqrt{s}$ such that \be\la{p23}
\Big(\int\limits_{Q(z_0,R_1)}|u|^4dz\Big)^\frac 12 <
R_1\sup\limits_{s/2<t<T} \|u(\cdot,t)\|_{4,\mathbb
R^2}<$$$$<R_1C(s/2,\|F\|_{4,Q_T},\|a\|_{2,\mathbb R^2})<\tau^4\ee
for any $z_0=(x_0,t_0)$ such that $t_0>s$. Obviously, $R_1$ depends
only  on $s$, $\|F\|_{q,Q_T}$, $\|F\|_{2,Q_T}$, and
$\|a\|_{2,\mathbb R^2}$ as, by interpolation,  $\|F\|_{4,Q_T}$ is
estimated by $\|F\|_{q,Q_T}$ and $\|F\|_{2,Q_T}$. Then we repeat the
proof of Theorem \ref{nrt1} replacing $Q$ with $Q(z_0,R_1)$ and
$Q_1$ with $Q(z_0,R_1/2)$ and establish $u\in
L^{\infty}(Q(z_0,R_1/2))$ with a uniform estimate with respect to
$z_0$ having $t_0>s$. Thus, the solution is H\"{o}lder continuous.

So, let us prove (\ref{p22}). To this end, we test the Navier-Stokes
by $|u|^2u$, as a result we have
$$\frac 14\pa_t \int\limits_{\mathbb R^2}|u|^4 dx
+\int\limits_{\mathbb R^2}\Big(|u|^2|\na u|^2 dx+2|u|^2|\na
|u||^2\Big) dx=$$
$$=\int\limits_{\mathbb R^2}\Big(2|u|pu\cdot  \na |u|+ F:(|u|^2\na u
+2|u|u\otimes \na |u|)\Big)dx. $$ After application of the Cauchy
inequality with a suitable weight, we find the following estimate
$$\pa_t\|u\|^4_{4,\mathbb R^2}\leq c \int\limits_{\mathbb R^2}\Big(p^2 |u|^2+
|F|^2|u|^2\Big)dx\leq$$$$\leq c\Big (\|p\|_{4,\mathbb
R^2}^2+\|F\|_{4,\mathbb R^2}^2\Big)\|u\|_{4,\mathbb R^2}^2.$$ It
remains to make use the pressure equation which, in the case of the
Cauchy problem, gives us \be\la{p24}\pa_t\|u(\cdot,t)\|^4_{4,\mathbb
R^2}\leq c\Big (\|u(\cdot,t)\|_{8,\mathbb
R^2}^4+\|F(\cdot,t)\|_{4,\mathbb
R^2}^2\Big)\|u(\cdot,t)\|_{4,\mathbb R^2}^2.\ee To evaluate the
right hand side of (\ref{p24}), we are going to use a particular
case of the generalized Ladyzhenskaya inequality
$$\|u(\cdot,t)\|_{8,\mathbb R^2}^2\leq c\|u(\cdot,t)\|_{4,\mathbb
R^2}\|\na u\|_{2,\mathbb R^2}.$$ The proof of the generalized
Ladyzhenskaya inequality is given in the Appendix. But then
(\ref{p24}) can be reduced to the form
$$\pa_t\|u(\cdot,t)\|^4_{4,\mathbb
R^2}\leq c\Big (\|u(\cdot,t)\|_{4,\mathbb R^2}^4\|\na u\|_{2,\mathbb
R^2}^2+\|F(\cdot,t)\|_{4,\mathbb R^2}^4+\|u(\cdot,t)\|_{4,\mathbb
R^2}^4\Big).$$ Multiplying the last inequality by a suitable cut-off
function in $t$, keeping in mind that our solution has the finite
energy bounded by $\|a\|^2_{2,\mathbb R^2}+\|F\|^2_{2,Q_T}$ from
above, and using Gronwall's lemma, we prove (\ref{p22}). Proposition
\ref{nrp4} is proved.

\setcounter{equation}{0}
\section{Appendix: Generalized Ladyzhenskaya Inequality}

The inequality \be\la{a1}\|u\|^2_{4,\mathbb R^2}\leq
\sqrt{2}\|u\|_{2,\mathbb R^2}\|\na u\|_{2,\mathbb R^2},\qquad
\forall u\in C^\infty_0(\mathbb R^2),\ee was used by Ladyzhenskaya
in \ci{La} to prove unique global solvability initial boundary value
problem for the Navier-Stokes equations in bounded domains of
$\mathbb R^2$. The generalized version of the Ladyzhenskaya
inequality is as follows: \be\la{a2} \|u\|^2_{2r,\mathbb R^2}\leq
\frac r{\sqrt{2}}\|u\|_{r,\mathbb R^2}\|\na u\|_{2,\mathbb
R^2},\qquad \forall u\in C^\infty_0(\mathbb R^2)\ee for $r\geq 2$.
The proof of (\ref{a2}) is essentially the same as the proof of
(\ref{a1}). The main ingredient of it the following identity
$$|u|^r(x_1,x_2)=r\int\limits^{x_1}_{-\infty}|u|^{r-1}(t,x_2)(|u|)_{,1}(t,x_2)dt
dx_2.$$ Using the same identity with respect to $x_2$, we find
$$\int\limits_{\mathbb R^2}|u|^{2r}dx_1dx_2\leq r^2\int\limits_{\mathbb R^2}
|u|^{r-1}|u_{,1}|dx_1dx_2\int\limits_{\mathbb
R^2}|u|^{r-1}|u_{,2}|dx_1dx_2\leq$$$$\leq \frac
{r^2}2\int\limits_{\mathbb R^2}
|u|^{2(r-1)}dx_1dx_2\int\limits_{\mathbb R^2} |\na u|^{2}dx_1dx_2.$$
By interpolation,
$$\int\limits_{\mathbb R^2}
|u|^{2(r-1)}dx_1dx_2\leq\Big(\int\limits_{\mathbb R^2}
|u|^{2r}dx_1dx_2\Big)^\frac {r-2}r\Big(\int\limits_{\mathbb R^2}
|u|^{r}dx_1dx_2\Big)^\frac 2r.$$ Now we deduce (\ref{a2}) from the
latter inequalities.\\
\vspace{1cm}

{\bf{Acknowledment}} P.C.'s research was partially sponsored by NSF grant DMS-0804380. G. S.'s research was partially supported by  the
RFFI grant 08-01-00372-a. P.C. gratefully acknowledges the hospitality of Oxford University's OxPDE Center.


\end{document}